\documentclass[12pt]{report}
\usepackage{amssymb,latexsym,epsfig,mathrsfs,graphpap,graphics,amsmath,amssymb}
\usepackage{amstext}
\usepackage{amsthm}
\usepackage[utf8]{inputenc}
\usepackage[english]{babel}
\usepackage{helvet}
\usepackage{courier}

\theoremstyle{Definition}

\theoremstyle{remark}

\numberwithin{equation}{section}
\begin{document}

\begin{quote}

{\bf\Large { On Quaternion Shearlet Transforms}}

\parindent=0mm \vspace{.4in}

  {\bf{Firdous A. Shah$^{\star}$ and Azhar Y. Tantary$^{\star}$ }}

 \parindent=0mm \vspace{.1in}
{\small \it $^{\star}$Department of  Mathematics,  University of Kashmir, South Campus, Anantnag-192101, Jammu and Kashmir, India. E-mail: $\text{fashah79@gmail.com}$;$\text{aytku92@gmail.com}$  }

{\small {\it Abstract:} In this paper, we introduce the notion of  quaternion shearlet transform- which is an extension of the ordinary shearlet transform. Firstly, we study the fundamental properties of quaternion shearlet transforms and then establish some basic results  including  Moyal's  and inversion formulae. Finally, we derive the associated Heisenberg's uncertainty inequality and the corresponding logarithmic version for quaternion shearlet transforms.

\parindent=0mm \vspace{.0in}

\parindent=0mm \vspace{.0in}
 }
\end{quote}

\parindent=0mm \vspace{.0in}
{\bf{1. Introduction}}

\parindent=0mm \vspace{.1in}
Since the inception of wavelets, their importance in the development of science and engineering, particularly in the areas of signal and image processing, is widely acknowledged {\cite M,\cite D}. Despite the fact that wavelet transform provides an optimal approximation for one dimensional data, it does not perform equally well in higher dimensions mainly due to the intrinsic isotropic nature. Several methodologies have been introduced in the recent literature to overcome these limitations of the traditional wavelet transform. In 2005, Kutyniok and Labate {\cite G} introduced the concept of shearlet transform for analyzing signals in higher dimensions. Unlike the wavelet transform, the shearlet transform is non-isotropic in nature and is ,therefore capable of  capturing the geometric features of signals including the detection of edges and surface discontinuities.

\parindent=8mm \vspace{.1in}
It is well known, shearlet theory is still in the developing phase and everyday many efforts are being made to extend this theory to a wider class of function spaces. On the  other hand,  considerable attention is also being  paid for the representation of signals in quaternion domain, so it is  quite interesting to extend the classical shearlet transforms to quaternion domains. The importance of the quaternion shearlet transform lies in the fact that it shall be capable of transforming a real 2D singal into a quaternion valued frequency domain while preserving the inherent features of the shearlet transform. Therefore, the main aim of this article is to introduce and investigate different properties of  quaternion shearlet transforms using the machinery of quaternion Fourier transforms.  Moreover, we derive the classical Heisenberg–Pauli–Weyl inequality and logarithmic version of this inequality {\cite F} for the quaternion shearlet transforms.

\parindent=0mm \vspace{.2in}
{\bf{2. Quaternion Algebra and Shearlet Transforms}}

\parindent=0mm \vspace{.1in}
The theory of quaternions was initiated by the Irish mathematician W.R Hamilton in 1842. The quaternion algebra provides an extension of the complex number system to a four dimensional non-commutative, associative algebra. The quaternion algebra is denoted by $\mathbb H$ and is given by
\begin{align*}
\mathbb H=\Big\{{\bf {h}}=a_{0}+i\,a_{1}+j\,a_{2}+k\,a_{3}\,:\,a_{0},a_{1},a_{2},a_{3}\in\mathbb R\Big\}
\end{align*}
where $i,j,k$ denote the three imaginary units, obeying the Hamilton's multiplication rules
\begin{align*}
ij=k=-ji,~jk=i=-kj\,,\,ki=j=-ik\,,\, {\text {and}}\,\, {i}^{2}={j}^{2}={k}^{2}={ijk}=-1.
\end{align*}
 Let ${\bf {h}}_{1}=a_{0}+i\,a_{1}+j\,a_{2}+k\,a_{3},{\bf {h}}_{2}=b_{0}+i\,b_{1}+ j\,b_{2}+k\,b_{3} {\text { be any two elements in }} \mathbb H,$ then the addition in $\mathbb H$ is defined componentwise and the multiplication in $\mathbb H$ is defined as
\begin{align*}
{\bf {h}}_{1}{\bf {h}}_{2}&=(a_{0}b_{0}-a_{1}b_{1}-a_{2}b_{2}-a_{3}b_{3})+i\,(a_{1}b_{0}+a_{0}b_{1}+a_{2}b_{3}-a_{3}b_{2})\\
&~~~+ j(a_{0}b_{2}+a_{2}b_{0}+a_{3}b_{1}-a_{1}b_{3})+k\,(a_{0}b_{3}+a_{3}b_{0}+a_{1}b_{2}-a_{2}b_{1}).
\end{align*}
Moreover, given any quaternion ${\bf {h}}=a_{0}+i\,a_{1}+j\,a_{2}+k\,a_{3},\,$ the conjugate is given by ${\overline{\bf {h}}}=a_{0}-i\,a_{1}-j\,a_{2}-k\,a_{3}$ and the norm is given by ${\|{\bf {h}}\|}_{\mathbb H}={\bf {h}}{\overline{\bf {h}}}=\sqrt{{a_{0}}^{2}+{a_{1}}^{2}+{a_{2}}^{2}+{a_{3}}^{2}}$. We note that any quaternion ${\bf {h}}=a_{0}+ i\,a_{1}+j\,a_{2}+k\,a_{3}$ can  be put in the form ${\bf {h}}=(a_{0}+i\,a_{1})+j\,(a_{2}-i\,a_{3})=u+j\,v$, where $u,v\in\mathbb C$. Also, ${\overline{\bf {h}}}={\overline u}-j\,v$, where $\overline{u}$ denotes the complex conjugate of $u$.
For any ${\bf {h}}_{1}=u_{1}+\,jv_{1},\, {\bf {h}}_{2}=u_{2}+j\,v_{2}$ in ${\mathbb H}$, the inner product is given by
\begin{align*}
{\langle {\bf {h}}_{1},{\bf {h}}_{2}\rangle}_{\mathbb H}&= {\bf {h}}_{1}\overline{{\bf {h}}_{2}}=(u_{1}{\overline u_{2}}+{\overline v_{1}}v_{2})+j(v_{1}{\overline u_{2}}-{\overline u_{1}}v_{2})
\end{align*}
We  see  that any function  $F$  from ${\mathbb R^{2}}$ to  ${\mathbb H}$ can be expressed as  $F(x)=f_{1}+j\,f_{2}$,  where $f_{1},f_{2}$ are both complex valued functions. Here we note that the following notations will be followed throughout this paper, $\check{F}(x)= F(-x)$ and $\tilde{F}(x)=\overline{f_{1}}(x)-j\check{f_{2}}(x).$

\pagestyle{myheadings}
\parindent=8mm \vspace{.1in}
The function space $L^2(\mathbb R^2,\mathbb H)$ consists of all measurable quaternion valued functions $F$ satisfying
\begin{align*}
 {\|F\|}_{2}=\left\{\int_{\mathbb R^{2}}\left(|f_{1}(x)|^{2}+|f_{2}(x)|^{2}\right)\,dx\right\}^{1/2}< \infty.
\end{align*}
Let $F=f_{1}+j\,f_{2},\,G=g_{1}+j\,g_{2}\in L^2(\mathbb R^2,\mathbb H)$,\, then the inner product in $L^2(\mathbb R^2,\mathbb H)$ is given by
\begin{align*}
{\langle F,G\rangle}_{2}&= \int_{\mathbb R^{2}}{\langle F,G\rangle}_{\mathbb H}\,dx\\&=\int_{\mathbb R^{2}}\bigg\{\left(f_{1}(x)\,\overline{g_{1}}(x)+\overline{f_{2}}(x)\,g_{2}(x)\right)+j\,\left(f_{2}(x)\,\overline{g_{1}}(x)
-\overline{f_{1}}(x)\,g_{2}(x)\right)\bigg\}\,dx.
\end{align*}
An easy computation shows that $L^2(\mathbb R^2,\mathbb H)$ equipped with  inner product defined above is a Hilbert space.

\parindent=0mm \vspace{.1in}
{\bf Definition 2.1.}  For any function $F\in L^1(\mathbb R^2,\mathbb H)\cap L^2(\mathbb R^2,\mathbb H)$, the   quaternion Fourier transform (QFT) is denoted by ${\mathscr F_{Q}}$ and is given by
\begin{align*}
{\mathscr F_{Q}}\big\{F(t)\big\}(\omega) =\hat F(\omega)= \int_{\mathbb R^2}e^{-2\pi {\bf i} t_{1} \cdot \omega_{1}} F(t)\,e^{-2\pi {\bf j} t_{2} \cdot \omega_{2}}\,dt,\tag{2.1}
\end{align*}
where $t=(t_{1},t_{2}),\, \omega=(\omega_{1},\omega_{2})$.
The  inverse quaternion Fourier transform corresponding to (2.1) is given by
\begin{align*}
F(t) &=\int_{\mathbb R^2}e^{2\pi {\bf i} t_{1} \cdot \omega_{1}} \hat F(\omega)\,e^{2\pi {\bf j} t_{2} \cdot \omega_{2}}\,d\omega.
\end{align*}

{\bf Definition 2.2.}  For any two complex valued functions $f,g\in L^2(\mathbb R^2)$, the convolution is given by
\begin{align*}
(f\ast g)(t)=\int_{\mathbb R^2}f(x)g(t-x)\, dx,\quad t\in \mathbb R^2.
\end{align*}
In analogy with the above definition, given any two quaternion valued functions\, $F ,G \in L^2(\mathbb R^2,\mathbb H)$  such that \,$F=f_{1}+jf_{2} ,\,\, G=g_{1}+j\,g_{2}$, the quaternion convolution is defined as
\begin{align*}
(F \star G)(t)= \left[\left(f_{1}\ast g_{1}\right)(t)-\left(\check{\overline {f}}_{2}\ast g_{2}\right)(t)\right]+ j \left[\left(\check{\overline {f_{1}}}\ast g_{2}\right)(t)+\left(f_{2}\ast g_{1}\right)(t)\right],t\in \mathbb R^2.
\end{align*}

\parindent=0mm \vspace{.0in}
 The following theorem provides some of the fundamental properties of QFT, the proof can be found in {\cite H   }.

\parindent=0mm \vspace{.1in}
{\bf Theorem 2.3.}  For any $F, G\in L^1(\mathbb R^2,\mathbb H)\cap L^2(\mathbb R^2,\mathbb H)$, the QFT has the following properties:

$(1)\,{\it Linearity:}\quad{\mathscr F_{Q}}\left(h_{1}F+h_{2}G\right)(\omega)=h_{1}{\mathscr F_{Q}}\left(F\right)(\omega)+
h_{2}\,{\mathscr F_{Q}}\left(G\right)(\omega),\quad h_{1},h_{2}\in {\mathbb H}.\\
(2)\,{\text {{\it Parseval formula:}}}\quad|{\mathscr F_{Q}}\left(F\right)(\omega)\|=\|F\|.\\
(3)\,{\text {{\it Convolution theorem:}}}\quad{\mathscr F_{Q}}\left[(F \star G)(t)\right](\omega)={\mathscr F_{Q}}\left[F\right](\omega)\,{\mathscr F_{Q}}\left[G\right](\omega)$.

\parindent=8mm \vspace{.1in}
  Before we proceed to establish the main results of this paper we shall briefly revise the continuous shearlet transform. Let $a\in\mathbb R^{+},s\in{\mathbb R}$ and $t\in\mathbb R^{2}$. Consider the square matrices  $A_{a}=\left(
  \begin{array}{cc}
   a & 0 \\
   0 & {\sqrt{a}} \\
   \end{array}
   \right)$   and $S_{s}=\left(
   \begin{array}{cc}
     1 & s \\
     0 & 1 \\
     \end{array}
      \right)$.
Then, the dilation, shear and translation operators are respectively given as:
\begin{align*}
\big(D_{A_{a}}\psi\big)(x)&=|A_{a}|^{-1/2}\,{\psi}({A_{a}^{-1}}x),\\
\big(D_{S_{s}}\psi\big)(x)&={\psi}({S_{s}^{-1}}x),\\
\big({T_{t}}\psi\big)(x)&={\psi}(x-t).
\end{align*}

\parindent=0mm \vspace{.0in}
Consequently, the family of shearlets is given by
\begin{align*}
 {\mathbb S}=\Big\{{\psi}_{a,s,t}(x)=|A_{a}|^{-1/2}\psi\left({A_{a}^{-1}}{S_{s}^{-1}}(x-t)\right):a\in\mathbb R^{+},s\in{\mathbb R}, t\in\mathbb R^{2}\Big\}.
\end{align*}

\parindent=0mm \vspace{.0in}
{\bf Definition 2.4.} Suppose that $\psi\in L^2(\mathbb R^2)$ satisfies the condition
\begin{align*}
C_{\psi}=\int_{\mathbb R}\int_{\mathbb R^{+}}{\left\|\hat{\psi}(\omega S_{s} A_{a})\right\|^{2}}\,\dfrac{da\,ds}{{a}^{3/2}}<\infty,\tag{2.2}
\end{align*}
 then we say that $\psi$ is an {\text{\it admissible shearlet}}.
\parindent=0mm \vspace{.1in}\\
{\bf Definition 2.5.} Let $\psi\in L^2(\mathbb R^2)$ be an {\text{\it admissible shearlet}}. For any $f\in L^2(\mathbb R^2)$, the continuous shearlet transform is defined as
\begin{align*}
S_{\psi}f(a,s,t)=\langle{f,\psi_{a,s,t}}\rangle=\int_{\mathbb R^{2}} f(x)\,\overline{\psi_{a,s,t}(x)}\,dx. \tag{2.3}
\end{align*}
{\bf Proposition 2.6.} Let $\psi\in L^2(\mathbb R^2)$ satisfies the admissibility condition (2.2). For any $f\in L^2(\mathbb R^2)$, the continuous shearlet transform can be represented as
\begin{align*}
{S_{\psi}}f(a,s,t)=\left(f\ast{\check{\overline{\psi}}}_{a,s,0}\right)(t),\quad t\in\mathbb R^2.
\end{align*}
{\it Proof.} We observe that
\begin{align*}
\psi_{a,s,t}(x)&=|A_{a}|^{-1/2}\,\psi\left({A_{a}}^{-1}{S_{s}}^{-1}(x-t)\right),\\
\psi_{a,s,0}(x)&=|A_{a}|^{-1/2}\,\psi\left({A_{a}}^{-1}{S_{s}}^{-1}x\right),\\
\check{\psi}_{a,s,0}(x)&=|A_{a}|^{-1/2}\,\psi\left({A_{a}}^{-1}{S_{s}}^{-1}(-x)\right).
\end{align*}
Furthermore, we have
\begin{align*}
\left(f\ast{\check{\overline{\psi}}}_{a,s,0}\right)(t)&=\int_{\mathbb R^{2}} f(x)\,{\check{\overline{\psi}}}_{a,s,0}(t-x)\,dx\\
&=|A_{a}|^{-1/2}\int_{\mathbb R^{2}} f(x)\overline{\psi\left({A_{a}}^{-1}{S_{s}}^{-1}(-(t-x))\right)}\\
&=|A_{a}|^{-1/2}\int_{\mathbb R^{2}} f(x)\,\overline{\psi\left({A_{a}}^{-1}{S_{s}}^{-1}(x-t)\right)}\\&=S_{\psi}f(a,s,t).
\end{align*}

\parindent=0mm \vspace{.0in}
{\bf Lemma 2.7.} Suppose that $f,\psi\in{L^2(\mathbb R^{2})}$ with $\psi$ satisfying (2.2), then we have
\begin{align*}
{\mathscr F}\big(S_{\psi}f(a,s,t)\big)\big(\omega\big)=|A_{a}|^{1/2}\,\hat{f}(\omega)\,\overline{\hat{\psi}(\omega S_{s}A_{a})}.
\end{align*}

\parindent=0mm \vspace{.0in}
{\it Proof.} By virtue of Plancheral theorem, we obtain
\begin{align*}
S_{\psi}f(a,s,t)&=\int_{\mathbb R^{2}}f(x)\,\overline{{\psi}_{a,s,t}(x)}\,dx\\
&=\int_{\mathbb R^{2}}\hat{f}(\omega)\,\overline{\hat{{\psi}}_{a,s,t}(\omega)}\,d{\omega}\\
&=\int_{\mathbb R^{2}}\hat{f}(\omega)\,\left(\int_{\mathbb R^{2}}\overline {|A_{a}|^{-1/2}\psi\left({A_{a}}^{-1} {S_{s}}^{-1}(x-t))\right)\,e^{-2\pi {\bf i}\omega \cdot x}\,dx}\right)\,d\omega\\
&=|A_{a}|^{1/2}\int_{\mathbb R^{2}}\hat{f}(\omega)\,\Big(\int_{\mathbb R^{2}}\overline{\psi(z)\,e^{-2\pi{\bf i}\omega\cdot(S_{s}A_{a}z+t)}\,dz}\Big)\,d\omega\\
&=|A_{a}|^{1/2}\int_{\mathbb R^{2}}\hat{f}(\omega)\,\left(\overline{\hat{\psi}(\omega S_{s}A_{a})\,e^{-2\pi {\bf i}\omega\cdot t}}\right)\,d\omega\\
&=|A_{a}|^{1/2}\int_{\mathbb R^{2}}\hat{f}(\omega)\,\overline{\hat{\psi}(\omega S_{s}A_{a})}\,e^{2\pi i \omega\cdot t}\,d\omega\\
&=|A_{a}|^{1/2}{\mathscr F}^{-1}\left(\hat{f}(\omega)\,\overline{\hat{\psi}(\omega S_{s}A_{a})}\right).
\end{align*}
Applying Fourier transform on both sides of above equation, we get
\begin{align*}
{\mathscr F}\left(S_{\psi}f(a,s,t)\right)\left(\omega\right)=|A_{a}|^{1/2}\hat{f}(\omega)\,\overline{\hat{\psi}(\omega S_{s}A_{a})}.
\end{align*}
This establishes the desired result.

\parindent=0mm \vspace{.1in}
{\it Observation:} A simple computation shows that, for $i=1,2$
\begin{align*}
{\mathscr F}\left[{\psi_{i}}_{a,s,t}\right](\omega)&=|A_{a}|^{1/2}\,{\mathscr F}\left[\psi_{i}\right](\omega\,S_{s}\,A_{a})\,e^{-2\pi i\omega\cdot t},\\
{\mathscr F}\left[{\check{\overline{\psi_{i}}}}_{a,s,t}\right](\omega)&=|A_{a}|^{1/2}\,{\mathscr F}\left[\check{\overline{\psi_{i}}}\right](\omega\,S_{s}\,A_{a})\,e^{-2\pi i\omega\cdot t}.
\end{align*}

\parindent=0mm \vspace{.0in}
{\bf{3. Quaternion Shearlet Transforms}}

\parindent=0mm \vspace{.1in}
In this Section, we shall introduce the notion of quaternion shearlet transforms and investigate some of their fundamental properties.

\parindent=0mm \vspace{.1in}
{\bf Definition 3.1.} If $\Psi\in L^2(\mathbb R^2,\mathbb H)$ satisfies the condition
\begin{align*}
C_{\Psi}=\int_{\mathbb R}\int_{\mathbb R}\left\|\hat{\Psi}({\omega}_{1},{\omega}_{2})\right\|_{\mathbb H}^{2}\,\dfrac{d{\omega}_{1}d{\omega}_{2}}{{\omega}_{1}^{2}}<\infty,\tag{3.1}
\end{align*}
then we say that $\Psi$ is admissible. The condition (3.1) is called the admissibility condition. Following {\cite S} and using the linearity property of the quaternion Fourier transforms, we can rewrite the admissibility condition (3.1) as
\begin{align*}
C_{\Psi}=\int_{\mathbb R}\int_{\mathbb R^{+}}{\left\|\hat{\Psi}(\omega S_{s} A_{a})\right\|_{\mathbb H}^{2}}\,\dfrac{da\,ds}{{a}^{3/2}}<\infty.\tag{3.2}
\end{align*}

\parindent=0mm \vspace{.1in}
For $\Psi\in L^2(\mathbb R^2,\mathbb H)$, we define the quaternion shearlet system as
\begin{align*}
\Psi_{a,s,t}(x)=|A_{a}|^{-1/2}\psi\left({A_{a}^{-1}}{S_{s}^{-1}}(x-t)\right), \quad a\in\mathbb R^{+},s\in\mathbb R,t\in\mathbb R^{2}.\tag{3.3}
\end{align*}
Corresponding to (3.3), we have the following definition of  quaternion shearlet transform.

\parindent=0mm \vspace{.1in}
{\bf Definition 3.2.} If $\Psi\in L^2(\mathbb R^2,\mathbb H)$ and $\Psi_{a,s,t}(x)$ is given by (3.3), then the integral transformation ${\mathcal S_{\Psi}}$ defined on $L^2(\mathbb R^2,\mathbb H)$ by
\begin{align*}
{\mathcal S_{\Psi}}F(a,s,t)&=\big\langle{F,\Psi_{a,s,t}}\big\rangle_{2}=\int_{\mathbb R^{2}} F(x)\overline{\Psi}_{a,s,t}(x)\,dx. \tag{3.4}.
\end{align*}
is called {\it quaternion shearlet transform} of $F(x)$.

\parindent=8mm \vspace{.1in}
By virtue of the  quaternion convolution, we can rewrite quaternion shearlet transform (3.4) as
\begin{align*}
{\mathcal S_{\Psi}}F(a,s,t)=\left(F\star{\check{\tilde{\Psi}}}_{a,s,0}\right)(t),\quad t\in\mathbb R^2 \tag{3.5}.
\end{align*}

\parindent=0mm \vspace{.0in}
We now prove a lemma which gives the decomposition of the quaternion shearlet transform in terms of the ordinary  shearlet transform.

\parindent=0mm \vspace{.1in}
{\bf Lemma 3.3.} If $F,\Psi\in L^2(\mathbb R^2,\mathbb H)$ , then the  quaternion shearlet transform (3.4) can be decomposed as
\begin{align*}
{\mathcal S}_{\Psi}F(a,s,t)=\left(S_{{\psi}_{1}}f_{1}+\overline{S_{\check{{\psi}_{2}}}\check{f_{2}}}\right)(a,s,t)+j\left(S_{{\psi}_{1}}f_{2}-\overline{ S_{\check{{\psi}_{2}}}\check{f_{1}}}\right)(a,s,t).
\end{align*}
{\it Proof.}   For any $t\in\mathbb R^2$, we observe that
\begin{align*}
{\mathcal S_{\Psi}}F(a,s,t)&=\left(F\star{\check{\tilde{\Psi}}}_{a,s,0}\right)(t)\\&=\left[\left(f_{1}+jf_{2}\right)\star\left(\check{\overline{{\psi_{1}}}}_{a,s,0}-j{\psi_{2}}_{a,s,0}\right)\right](t)\\&=
\left(f_{1}\ast\check{{\overline{\psi_{1}}}}_{a,s,0}\right)(t)+\left(\check{\overline{f_{2}}}\ast{\psi_{2}}_{a,s,0}\right)(t)
+j\left[\left(f_{2}\ast\check{{\overline{\psi_{1}}}}_{a,s,0}\right)-\left(\check{\overline{f_{1}}}\ast{\psi_{2}}_{a,s,0}\right)\right](t)\tag{3.6}
\end{align*}
For $k=1,2$, we have
\begin{align*}
\left(\check{\overline{f_{k}}}\ast{\psi_{2}}_{a,s,0}\right)(t)&=\int_{\mathbb R^{2}}\check{\overline{f_{k}}}(x)\,{\psi_{2}}_{a,s,0}(t-x)\,dx\\&=\int_{\mathbb R^{2}}\check{\overline{f_{k}}}(x)\,\check{{\psi}_{2}}_{a,s,0}(x-t)\,dx\\&=\int_{\mathbb R^{2}}\overline{\check{f_{k}}(x)\,\overline{{\check{{\psi}_{2}}}}_{a,s,0}(x-t)}\,dx\\&=\overline{ S_{\check{{\psi}_{2}}}\check{f_{k}}}(a,s,t).
\end{align*}
Plugging the above estimates in (3.6) and then using Proposition 2.6, we obtain
\begin{align*}
{\mathcal S}_{\Psi}F(a,s,t)=\left(S_{{\psi}_{1}}f_{1}+\overline{S_{\check{{\psi}_{2}}}\check{f_{2}}}\right)(a,s,t)+j\left(S_{{\psi}_{1}}f_{2}-\overline{ S_{\check{{\psi}_{2}}}\check{f_{1}}}\right)(a,s,t).
\end{align*}

\parindent=0mm \vspace{.0in}
In the following theorem, we assemble some of the basic properties of the  quaternion shearlet transform.

\parindent=0mm \vspace{.1in}
{\bf Theorem 3.4.} {\it For any  $F\,,G\in L^2(\mathbb R^2,\mathbb H)$, the  quaternion shearlet transform (3.4) satisfies the following properties:}

\parindent=0mm \vspace{.1in}
$(i)~~{\text {\it Linearity:}}\quad{\mathcal S}_{\Psi}\left(h_{1}F+h_{2}G\right)(a,s,t)=h_{1}{\mathcal S}_{\Psi}F(a,s,t)+h_{2}{\mathcal S}_{\Psi}(a,s,t),\quad  h_{1}, h_{2}\in{\mathbb H}$

\parindent=0mm \vspace{.1in}
$(ii)~~{\text {\it Anti-linearity:}}\quad{\mathcal S}_{h_{1}\Psi+h_{2}\Phi}F(a,s,t)=\overline{h_{1}}\,{\mathcal S}_{\Psi}F(a,s,t)+\overline{h_{2}}\,{\mathcal S}_{\Phi}(a,s,t)$

\parindent=0mm \vspace{.1in}
$(iii)~{\text {\it Translation:}}\quad{\mathcal S}_{\Psi}\left(T_{t^\prime}F\right)(a,s,t)=\left[S_{{\psi}_{1}}f_{1}\right](a,s,t- t^{\prime})+\overline{\left[S_{\check{{\psi}_{2}}}\check{f_{2}}\right]}(a,s,t+t^\prime)$

$\qquad\qquad\qquad\qquad\qquad\qquad\qquad+j\left(\left[ S_{{\psi}_{1}}f_{2}\right](a,s,t- t^{\prime})-\overline{\left[S_{\check{{\psi}_{2}}}\check{f_{1}}\right]}(a,s,t+ t^{\prime})\right), t^{\prime}\in{\mathbb R^{2}}$

\parindent=0mm \vspace{.0in}
$(iv)\,{\text {\it Scaling:}}\quad{\mathcal S}_{\Psi}\left(F(\lambda\, x)\right)(a,s,t)=\dfrac{1}{\lambda}\left[{\mathcal S}_{\Psi}F(x)\right](a,s/{\lambda},{\lambda\,t}),\quad \lambda\in{\mathbb R}$.

\parindent=0mm \vspace{.1in}
{\it Proof.} (i). For $ h_{1}, h_{2}\in{\mathbb H}$, we have
\begin{align*}
{\mathcal S}_{\Psi}\left(h_{1}F+h_{2}G\right)(a,s,t)&= \left({\left(h_{1}F+h_{2}G\right)}\star{\check{\tilde{\Psi}}}_{a,s,0}\right)(t)\\&=\int_{\mathbb R^{2}}\left\{{\left(h_{1}F+h_{2}G\right)(x)}{\check{\tilde{\Psi}}}_{a,s,0}(t-x)\right\}\,dx\\&=\int_{\mathbb R^{2}}h_{1}F(x)\,{\check{\tilde{\Psi}}}_{a,s,0}(t-x)\,dx+\int_{\mathbb R^{2}}h_{2}\,G(x)\,{\check{\tilde{\Psi}}}_{a,s,0}(t-x)\,dx\\&=h_{1}\left(F\star{\check{\tilde{\Psi}}}_{a,s,0}\right)(t)
+h_{2}\left(G\star{\check{\tilde{\Psi}}}_{a,s,0}\right)(t)\\&=h_{1}{\mathcal S}_{\Psi}F(a,s,t)+h_{2}{\mathcal S}_{\Psi}(a,s,t)
\end{align*}
(ii). For $ h_{1}, h_{2}\in{\mathbb H}$, we have
\begin{align*}
{\mathcal S}_{h_{1}\Psi+h_{2}\Phi}F(a,s,t)&= \left(F(x)\star\check{\widetilde{\left\{{h_{1}{\Psi}_{(a,s,0)}+h_{2}{\Phi}_{a,s,0}}\right\}}}\right)(t)\\&=\int_{\mathbb R^{2}}F(x)\,\check{\widetilde{\left\{{h_{1}\,{\Psi}_{(a,s,0)}+h_{2}\,{\Phi}_{a,s,0}}\right\}}}(t-x)\,dx\\&=\int_{\mathbb
R^{2}}\overline{h_{1}}\,F(x)\,{\check{\tilde{\Psi}}}_{a,s,0}(t-x)\,dx+\int_{\mathbb R^{2}}\overline{h_{2}}\,F(x)\,{\check{\tilde{\Phi}}}_{a,s,0}(t-x)\,dx\\&=\overline{h_{1}}\left(F\star{\check{\tilde{\Psi}}}_{a,s,0}\right)(t)
+\overline{h_{2}}\left(F\star{\check{\tilde{\Phi}}}_{a,s,0}\right)(t)\\&=\overline{h_{1}}\,{\mathcal S}_{\Psi}F(a,s,t)+\overline{h_{2}}\,{\mathcal S}_{\Phi}(a,s,t)
\end{align*}
(iii). For $t^{\prime}\in{\mathbb R^{2}}$, we have
\begin{align*}
&{\mathcal S}_{\Psi}\left(T_{ t^{\prime}}F\right)(a,s,t)\\
&\quad=\left( S_{{\psi}_{1}}(T_{ t^{\prime}}f_{1})+\overline{ S_{\check{{\psi}_{2}}}\check{(T_{ t^{\prime}}f_{2})}}\right)(a,s,t)+j\left({S_{{\psi}_{1}}}({T_{ t^{\prime}}}f_{2})-\overline{ S_{\check{{\psi}_{2}}}\check{{(T_{ t^{\prime}}f_{1})}}}\right)(a,s,t).\tag{3.7}
\end{align*}
Moreover, for $k=1,2$, we have
\begin{align*}
S_{{\psi}_{1}}T_{ t^{\prime}}f_{k}(a,s,t)&=\int_{\mathbb R^{2}}f_{k}(x- t^{\prime}){\overline{{\psi_{1}}_{a,s,t}(x)}}\,dx\\
&=|A_{a}|^{-1/2}\int_{\mathbb R^{2}}f_{k}(x- t^{\prime}){\overline{\psi_{1}({A_{a}}^{-1}{S_{s}}^{-1}(x-t))}}\,dx\\
&=|A_{a}|^{-1/2}\int_{\mathbb R^{2}}f_{k}(z){\overline{\psi_{1}({A_{a}}^{-1}{S_{s}}^{-1}((z-(t- t^{\prime}))}}\,dz\\&=S_{{\psi}_{1}}f_{k}\left(a,s,t- t^{\prime}\right)
\end{align*}
and
\begin{align*}
S_{\check{{\psi}_{2}}}\check{(T_{ t^{\prime}}f_{k})}(a,s,t)&=\int_{\mathbb R^{2}}\check{(T_{ t^{\prime}}f_{k})}(x)\,{\overline{\check{{\psi_{2}}}_{a,s,t}(x)}}\,dx\\
&=|A_{a}|^{-1/2}\int_{\mathbb R^{2}}(T_{ t^{\prime}}f_{k})(-x)\,{\overline{\check{\psi_{2}}({A_{a}}^{-1}{S_{s}}^{-1}(x-t))}}\,dx\\
&=|A_{a}|^{-1/2}\int_{\mathbb R^{2}}(T_{ t^{\prime}}f_{k})(y)\,{\overline{\check{\psi_{2}}({A_{a}}^{-1}{S_{s}}^{-1}(-y-t))}}\,(-dy)\\
&=|A_{a}|^{-1/2}\int_{\mathbb R^{2}}{f_{k}}(y- t^{\prime})\,{\overline{\check{\psi_{2}}({A_{a}}^{-1}{S_{s}}^{-1}(-y-t))}}\,(-dy)\\
&=|A_{a}|^{-1/2}
\int_{\mathbb R^{2}}{f_{k}}(-z)\,{\overline{\check{\psi_{2}}({A_{a}}^{-1}{S_{s}}^{-1}(z- t^{\prime}-t))}}\,\,dz\\
&=|A_{a}|^{-1/2}\int_{\mathbb R^{2}}\check{f_{k}}(z)\,{\overline{\check{\psi_{2}}({A_{a}}^{-1}{S_{s}}^{-1}(z-( t^{\prime}+t))}}\,\,dz\\
&= S_{\check{{\psi}_{2}}}\check{f_{k}}\left(a,s,t+ t^{\prime}\right).
\end{align*}
Using the above two estimates in (3.7), we get
\begin{align*}
{\mathcal S}_{\Psi}\left(T_{ t^{\prime}}F\right)(a,s,t)=\left[S_{{\psi}_{1}}f_{1}\right](a,s,t- t^{\prime})+\overline{\left[ S_{\check{{\psi}_{2}}}\check{f_{2}}\right]}(a,s,t+ t^{\prime})\\+j\left(\left[S_{{\psi}_{1}}f_{2}\right](a,s,t- t^{\prime})-\overline{\left[ S_{\check{{\psi}_{2}}}\check{f_{1}}\right]}\left(a,s,t+ t^{\prime}\right)\right)
\end{align*}

(iv). For $\lambda\in{\mathbb R}$, we have
\begin{align*}
&{\mathcal S}_{\Psi}\left(F(\lambda\, x)\right)(a,s,t)\\
&=\left(S_{{\psi}_{1}}f_{1}(\lambda\, x)+\overline{ S_{\check{{\psi}_{2}}}\check{f_{2}}(\lambda\, x)}\right)(a,s,t)+j\left(S_{{\psi}_{1}}{f_{2}}(\lambda\, x)-\overline{S_{\check{{\psi}_{2}}}\check{f_{1}}(\lambda\, x)}\right)(a,s,t).\tag{3.8}
\end{align*}
For $k=\,1,2$, we observe that
\begin{align*}
\left[{\mathcal S}_{{\psi}_{1}}f_{k}(\lambda\,x)\right](a,s,t)&=\int_{\mathbb R^{2}}f_{k}(\lambda\,x){\overline{{\psi_{1}}_{a,s,t}(x)}}\,dx\\
&=|A_{a}|^{-1/2}\int_{\mathbb R^{2}}f_{k}(\lambda\,x){\overline{\psi_{1}({A_{a}}^{-1}{S_{s}}^{-1}(x-t))}}\,dx\\
&=|A_{a}|^{-1/2}\dfrac{1}{\lambda}\int_{\mathbb R^{2}}f_{k}(z){\overline{\psi_{1}({A_{a}}^{-1}{S_{s}}^{-1}(z/{\lambda}-t))}}\,dz\\
&=|A_{a}|^{-1/2}\dfrac{1}{\lambda}\int_{\mathbb R^{2}}f_{k}(z){\overline{\psi_{1}({A_{a}}^{-1}{S_{s/\lambda}}^{-1}(z-\lambda t))}}\,dz\\&=\dfrac{1}{\lambda}\left[{\mathcal S}_{{\psi}_{1}}f_{k}(x)\right](a,s/{\lambda},{\lambda\,t}).
\end{align*}
A similar computation will show that
\begin{align*}
  \left[{\mathcal S}_{\check{{\psi}_{2}}}\check{f_{i}}{(\lambda\,x)}\right](a,s,t)=\dfrac{1}{\lambda}\left[{\mathcal S}_{\check{{\psi}_{2}}}\check{f_{i}}(x)\right](a,s/{\lambda},{\lambda\,t}),\,\qquad i=1,\,2.
\end{align*}
We note that   $S_{s/\lambda}=\left(
\begin{array}{cc}
1/{\lambda} & s/{\lambda} \\
0 & 1/{\lambda} \\
 \end{array}
 \right).$
Using the above two estimates in (3.8), we observe that
\begin{align*}
{\mathcal S}_{\Psi}\left(F(\lambda\, x)\right)(a,s,t)=\dfrac{1}{\lambda}\left[{\mathcal S}_{\Psi}F(x)\right](a,s/{\lambda},{\lambda\,t}).
\end{align*}

{\bf Lemma 3.5.} {\it Let $\Psi\in L^2(\mathbb R^2,\mathbb H)$ be an admissible shearlet. Then, for any $F\in L^2(\mathbb R^2,\mathbb H)$, we have}
\begin{align*}
{\mathscr F_{Q}}\left[{\mathcal S_{\Psi}}F(a,s,t)\right](\omega)={|A_{a}|}^{1/2}\,\hat{F}(\omega)\,\hat{\Psi}(\omega S_{s} A_{a}).
\end{align*}
{\it Proof.} By invoking the convolution theorem for QFT, we have
\begin{align*}
{\mathscr F_{Q}}\left[{\mathcal S_{\Psi}}F(a,s,t)\right](\omega)&={\mathscr F_{Q}}\left[(F \star \check{\tilde{{\Psi}}}_{a,s,0})(t)\right](\omega)\\&={\mathscr F_{Q}}\left[F\right](\omega)\,\,{\mathscr F_{Q}}\left[\check{\tilde{{\Psi}}}_{a,s,0}\right](\omega)\\&={\mathscr F_{Q}}\left[F\right](\omega)\,\,{\mathscr F_{Q}}\left[\check{\overline{{\psi_{1}}}}_{a,s,0}-j{\psi_{2}}_{a,s,0}\right](\omega) \\&={\mathscr F_{Q}}\left[F\right](\omega)\,\left[{\mathscr F_{Q}}[\check{\overline{{\psi_{1}}}}_{a,s,0}](\omega)-j\,{\mathscr F_{Q}}[{{\psi}_{2}}_{a,s,0}](\omega)\right]\\&={\mathscr F_{Q}}\left[F\right](\omega)\,\left\{{|A_{a}|}^{1/2}\,{\mathscr F_{Q}}\left[{\check{\overline{\psi}_{1}}}\right](\omega S_{s} A_{a})-j\,{|A_{a}|}^{1/2}\,{\mathscr F_{Q}}\left[{\psi}_{2}\right](\omega S_{s} A_{a})\right\}\\&={|A_{a}|}^{1/2}\,{\mathscr F_{Q}}\left[F\right](\omega)\,{\mathscr F_{Q}}\left[{\check{\overline{\psi}_{1}}}-j\,{\psi}_{2}\right](\omega S_{s} A_{a})\\&={|A_{a}|}^{1/2}\,{\mathscr F_{Q}}\left[F\right](\omega)\,{\mathscr F_{Q}}\left[\check{\tilde{{\Psi}}}\right](\omega S_{s} A_{a})\\&={|A_{a}|}^{1/2}\,\hat{F}(\omega)\,\overline{\hat\Psi(\omega S_{s} A_{a})},
\end{align*}
which evidently completes the proof of the Lemma.

\parindent=0mm \vspace{.1in}

{\bf Theorem 3.6} {(\bf Moyal's Principle).} {\it Suppose $\Psi\in L^2(\mathbb R^2,\mathbb H)$ satisfies the admissibility condition (3.1). Then, for any $F\,,G\in L^2(\mathbb R^2,\mathbb H)$, we have}
\begin{align*}
\int_{\mathbb R^2}\int_{\mathbb R}\int_{\mathbb R^{+}}{\mathcal S_{\Psi}}F(a,s,t)\,\overline{{{\mathcal S_{\Psi}}G(a,s,t)}}\dfrac{da\,ds\,dt}{a^{3}}=C_{\Psi}\,\big\langle F,G \big\rangle_{2}.
\end{align*}
{\it Proof.}  Consider two functions $F\,,G\in L^2(\mathbb R^2,\mathbb H)$ and suppose that  $\Psi\in L^2(\mathbb R^2,\mathbb H)$ is admissible. Then, by  Fubini theorem and   Plancheral theorem for QFT, we obtain
\begin{align*}
&\int_{\mathbb R^2}\int_{\mathbb R}\int_{\mathbb R^{+}}{\mathcal S_{\Psi}}F(a,s,t)\,\overline{{{\mathcal S_{\Psi}}G(a,s,t)}}\dfrac{da\,ds\,dt}{a^{3}}\\
&= \int_{\mathbb R}\int_{\mathbb R^{+}}\left\{\int_{\mathbb R^2}{\mathscr F_{Q}}\left[{\mathcal S_{\Psi}}F(a,s,t)\right](\omega)\,\overline{{\mathscr F_{Q}}\left[{\mathcal S_{\Psi}}F(a,s,t)\right](\omega)}\,d{\omega} \right\}\dfrac{da\,ds}{a^{3}}\\
&= \int_{\mathbb R}\int_{\mathbb R^{+}}\left\{\int_{\mathbb R^2}{\mid A_{a}\mid}^{1/2}\,\hat{F}(\omega)\,\overline{\hat\Psi(\omega S_{s} A_{a})}\,\,\,\overline{{\mid A_{a}\mid}^{1/2}\,\hat{G}(\omega)\,}\hat\Psi(\omega S_{s} A_{a})\,d{\omega}\right\}\dfrac{da\,ds}{a^{3}}\\
&= \int_{\mathbb R}\int_{\mathbb R^{+}}\left\{\int_{\mathbb R^2}\hat{F}(\omega)\overline{\hat{G}(\omega)}\,\,{\|\hat{\Psi}(\omega S_{s} A_{a})\|_{\mathbb H}^{2}} \right\}\dfrac{da\,ds}{a^{3}}\\&=\int_{\mathbb R^2}\left\{\int_{\mathbb R}\int_{\mathbb R^{+}}\dfrac{{\|\hat{\Psi}(\omega S_{s} A_{a})\|_{\mathbb H}^{2}}}{a^{3/2}}\,da\,ds \right\}\hat{F}(\omega)\overline{\hat{G}(\omega)}\,d{\omega}\\
&= C_{\Psi}\,\int_{\mathbb R^2}\hat{F}(\omega)\overline{\hat{G}(\omega)}\,d{\omega}\\&= C_{\Psi}\,\int_{\mathbb R^2}F(t)\overline{\hat{G}(t)}\,d{t}\\&=C_{\Psi}\,{\langle F,G \rangle}_{2},
\end{align*}
where $C_{\Psi}$ is given by (3.1). This completes the proof of the theorem.

\parindent=0mm \vspace{.1in}
{\bf Corollary 3.7.} {\it If $\Psi\in L^2(\mathbb R^2,\mathbb H)$ satisfies the admissibility condition (3.1). Then, for any any $F\in L^2(\mathbb R^2,\mathbb H)$, we have}
\begin{align*}
\int_{\mathbb R^2}\int_{\mathbb R}\int_{\mathbb R^{+}}\left\|\mathcal S_{\Psi}F(a,s,t)\right\|_{\mathbb H}^{2}\dfrac{da\,ds\,dt}{a^{3}}=C_{\Psi}\, \big\|F\big\|_{2}.\tag{3.9}
\end{align*}

\parindent=0mm \vspace{.0in}
From  (3.9),  it is clear to observe that the quaternion shearlet transform is a bounded linear operator from $L^2(\mathbb R^2,\mathbb H)$ to $L^2(\mathbb R^2\times\mathbb R\times\mathbb R^{+},\mathbb H)$.

\parindent=0mm \vspace{.1in}
{\bf Theorem 3.8 ({\bf Inversion Formula}).} {\it Let $F\in L^2(\mathbb R^2,\mathbb H)$ and ${\mathcal S_{\Psi}}F(a,s,t)$ denotes the quaternion shearlet transform of $F$ with respect to $\Psi\in L^2(\mathbb R^2,\mathbb H)$, then}
\begin{align*}
F(x)=\dfrac{1}{C_{\Psi}}\int_{\mathbb R}\int_{\mathbb R^{+}}\big({\mathcal S_{\Psi}}F(a,s,t)\star{\Psi}_{a,s,0}\big)(x)\,\dfrac{da\,ds}{a^{3}},\quad x\in{\mathbb R^2}.\tag{3.10}
\end{align*}

{\it Proof.} Using  convolution theorem for QFT and Lemma 2.7, we have
\begin{align*}
{\mathscr F_{Q}}\left[{\mathcal S_{\Psi}}F(a,s,t)\right](\omega)&={\mathscr F_{Q}}\left[(F \star \check{\tilde{{\Psi}}}_{a,s,0})(t)\right](\omega)\\&=|A_{a}|^{1/2}\,\hat {F}(\omega)\,\overline{\hat{\Psi}(\omega S_{s} A_{a})}.
\end{align*}
Moreover, we have
\begin{align*}
{\mathscr F_{Q}}\left[\left(F \star \check{\tilde{{\Psi}}}_{a,s,0}\right)\star\Psi_{a,s,0}\right](\omega)&={\mathscr F_{Q}}\left[F \star \check{\tilde{{\Psi}}}_{a,s,0}\right](\omega){\mathscr F_{Q}}\left[\Psi_{a,s,0}\right](\omega)\\&=|A_{a}|^{1/2}\,\hat {F}(\omega)\,\overline{\hat{\Psi}(\omega S_{s} A_{a})}\,|A_{a}|^{1/2}\,\hat{\Psi}(\omega S_{s} A_{a})\\
&=a^{3/2}\,\hat {F}(\omega)\,{\left\|\hat{\Psi}(\omega S_{s} A_{a})\right\|_{\mathbb H}^{2}}.
\end{align*}
 Integrating the above equation on both sides with respect to the measure $\dfrac{da\,ds}{a^{3}}$, we obtain
\begin{align*}
\int_{\mathbb R}\int_{\mathbb R^{+}}{\mathscr F_{Q}}\left[\left(F \star \check{\tilde{{\Psi}}}_{a,s,0}\right)\star\Psi_{a,s,0}\right](\omega)\dfrac{da\,ds}{a^{3}}&=\int_{\mathbb R}\int_{\mathbb R^{+}}a^{3/2}\,\hat {F}(\omega)\,{\|\hat{\Psi}(\omega S_{s} A_{a})\|_{\mathbb H}^{2}}\,\dfrac{da\,ds}{a^{3}}\\&=\hat {F}(\omega)\,C_{\Psi}.
\end{align*}
By virtue of  Fubini theorem and the definition of  inverse QFT, we have
\begin{align*}
F(x)&=\int_{\mathbb R^{2}}\hat{F}(\omega)\,e^{2\pi i \omega\cdot x}\,d\omega\\&=\dfrac{1}{C_\Psi}\int_{\mathbb R^{2}}\left\{ \int_{\mathbb R}\int_{\mathbb R^{+}}{\mathscr F_{Q}}\left[\left(F \star \check{\tilde{{\Psi}}}_{a,s,0}\right)\star\Psi_{a,s,0}\right](\omega)\dfrac{da\,ds}{a^{3}}\right\}\,e^{2\pi i \omega\cdot x}\,d\omega\\&=
\int_{\mathbb R}\int_{\mathbb R^{+}}\left\{\int_{\mathbb R^{2}}{\mathscr F_{Q}}\left[\left(F \star \check{\tilde{{\Psi}}}_{a,s,0}\right)\star\Psi_{a,s,0}\right](\omega)\,e^{2\pi i\omega \cdot x}\,d\omega\right\}\dfrac{da\,ds}{a^{3}}\\&=\int_{\mathbb R}\int_{\mathbb R^{+}} \left[\left(F \star \check{\tilde{{\Psi}}}_{a,s,0}\right)\star\Psi_{a,s,0}\right](x)\,\dfrac{da\,ds}{a^{3}}\\&=\dfrac{1}{C_{\Psi}}\int_{\mathbb R}\int_{\mathbb R^{+}}\left({\mathcal S_{\Psi}}F(a,s,t)\star{\Psi}_{a,s,0}\right)(x)\,\dfrac{da\,ds}{a^{3}},\quad x\in{\mathbb R^2},
\end{align*}
which evidently completes the proof.

\parindent=0mm \vspace{.1in}

{\bf 4.  Uncertainty Principles for  Quaternion Shearlet Transforms}

\parindent=0mm \vspace{.1in}
In this Section, we shall establish an analogue of the well-known Heisenberg's uncertainty inequality and the corresponding logarithmic version for the quaternion shearlet transform as defined in (3.4).

\parindent=0mm \vspace{.1in}
{\bf Theorem 4.1.} {\it Suppose $\Psi\in L^2(\mathbb R^2,\mathbb H)$  is an admissible shearlet, then for any non-zero function $F\in L^2(\mathbb R^2,\mathbb H)$, we have}
\begin{align*}
{\left(\int_{\mathbb R^{2}} \int_{\mathbb R}\int_{\mathbb R^{+}}\|t\|^{2}\,\|{\mathcal S_{\Psi}}F(a,s,t)\|_{\mathbb H}^{2}\dfrac{da\,ds\,dt}{a^{3}}\right)}^{1/2}\,\left(\int_{\mathbb R^{2}}\|\omega\|^{2}\,{\|{\hat F}(\omega)\|_{\mathbb H}^{2}}\right)^{1/2}\geq\,\dfrac{\sqrt{C_{\Psi}}}{2\pi}{\|F\|}_{2}^{2}.
\end{align*}

\parindent=0mm \vspace{.1in}
{\it Proof.}  It is well known that for any quaternion valued function $F\in L^2(\mathbb R^2,\mathbb H)$, the Heisenberg Paul-weyl inequality  {\cite C} is given by
\begin{align*}
{\left(\int_{\mathbb R^{2}}\|t\|^{2}\,\|F(t)\|_{\mathbb H}^{2}\,dt\right)}^{1/2}\,\left(\int_{\mathbb R^{2}}\|\omega\|^{2}\,{\|{\hat F}(\omega)\|_{\mathbb H}^{2}}\right)^{1/2}\geq\,\dfrac{1}{2\pi}\,\int_{\mathbb R^{2}}\|F(t)\|_{\mathbb H}^{2}\,dt.\tag{4.1}
\end{align*}
As ${\mathcal S_{\Psi}}F(a,s,t)\in L^2(\mathbb R^2,\mathbb H)$ whenever $F\in L^2(\mathbb R^2,\mathbb H)$, so we can replace $F$ by  ${\mathcal S_{\Psi}}F(a,s,t)$ in the inequality (4.1) to obtain
\begin{align*}
{\left(\int_{\mathbb R^{2}}\|t\|^{2}\,\|{\mathcal S_{\Psi}}F(a,s,t)\|_{\mathbb H}^{2}\,dt\right)}^{1/2}\,\left(\int_{\mathbb R^{2}}\|\omega\|^{2}\,{\left\|\mathscr F_{Q}\left[\mathcal S_{\Psi}F(a,s,t)\right]\right\|_{\mathbb H}^{2}}\right)^{1/2}\\\geq\dfrac{1}{2\pi}\,\int_{\mathbb R^{2}}\|{\mathcal S_{\Psi}}F(a,s,t)\|_{\mathbb H}^{2}\,dt.
\end{align*}
We integrate the above inequality with respect to measure $da\,ds/a^{3}$,\, so that
\begin{align*}
\int_{\mathbb R}\int_{\mathbb R^{+}}\left\{{\left(\int_{\mathbb R^{2}}\|t\|^{2}\,\|{\mathcal S_{\Psi}}F(a,s,t)\|_{\mathbb H}^{2}\,dt\right)}^{1/2}\,\left(\int_{\mathbb R^{2}}\|\omega\|^{2}\,{\left\|\mathscr F_{Q}\left[\mathcal S_{\Psi}F(a,s,t)\right]\right\|_{\mathbb H}^{2}}\right)^{1/2}\right\}\dfrac{ds\,dt}{a^{3}}\\\geq\,\dfrac{1}{2\pi}\,\int_{\mathbb R}\int_{\mathbb R^{+}}\int_{\mathbb R^{2}}\|{\mathcal S_{\Psi}}F(a,s,t)\|_{\mathbb H}^{2}\,\dfrac{da\,ds\,dt}{a^{3}}.
\end{align*}
Using  Schwartz inequality followed by  Fubini theorem and (3.9), we obatin
\begin{align*}
&{\left(\int_{\mathbb R^{2}} \int_{\mathbb R}\int_{\mathbb R^{+}}\|t\|^{2}\,\|{\mathcal S_{\Psi}}F(a,s,t)\|_{\mathbb H}^{2}\,\dfrac{da\,ds\,dt}{a^{3}}\right)}^{1/2} \\
&\qquad\quad{\left(\int_{\mathbb R^{2}} \int_{\mathbb R}\int_{\mathbb R^{+}}\|\omega\|^{2}\,\left\|\mathscr F_{Q}\left[\mathcal S_{\Psi}F(a,s,t)\right]\right\|_{\mathbb H}^{2}\,\dfrac{da\,ds\,dt}{a^{3}}\right)}^{1/2} \geq\,\dfrac{1}{2\pi}\,C_{\Psi}\,\|F\|_{2}^{2}\tag{4.2}.
\end{align*}
Invoking  Lemma 3.5 and using (3.1), we observe that
\begin{align*}
&\int_{\mathbb R^{2}} \int_{\mathbb R}\int_{\mathbb R^{+}}\|\omega\|^{2}\,\left\|\mathscr F_{Q}\left[\mathcal S_{\Psi}F(a,s,t)\right]\right\|_{\mathbb H}^{2}\,\dfrac{da\,ds\,dt}{a^{3}}\\
&=\int_{\mathbb R^{2}} \int_{\mathbb R}\int_{\mathbb R^{+}}\|\omega\|^{2}\left\|\hat{F}(\omega)\right\|_{\mathbb H}^{2}\,\|\hat{\Psi}(\omega S_{s} A_{a})\|_{\mathbb H}^{2}\,\dfrac{da\,ds\,d\omega}{a^{3/2}}\\
&=\int_{\mathbb R^{2}}\Big\{\int_{\mathbb R}\int_{\mathbb R^{+}}\dfrac{\left\|\hat{\Psi}(\omega S_{s} A_{a})\right\|_{\mathbb H}^{2}}{a^{3/2}}\,da\,ds\Big\}\|\omega\|^{2}\,\|\hat{F}(\omega)\|_{\mathbb H}^{2}\,d\omega\\
&=C_{\Psi}\int_{\mathbb R^{2}}\|\omega\|^{2}\,\|\hat{F}(\omega)\|_{\mathbb H}^{2}\,d\omega.
\end{align*}
Using the above estimate in (4.2), we obtain the desired result 
\begin{align*}
{\left(\int_{\mathbb R^{2}} \int_{\mathbb R}\int_{\mathbb R^{+}}\|t\|^{2}\,\|{\mathcal S_{\Psi}}F(a,s,t)\|_{\mathbb H}^{2}\,\dfrac{da\,ds\,dt}{a^{3}}\right)}^{1/2}\,\left(\int_{\mathbb R^{2}}\|\omega\|^{2}\,{\|{\hat F}(\omega)\|_{\mathbb H}^{2}}\right)^{1/2}\geq\,\dfrac{\sqrt{C_{\Psi}}}{2\pi}\|F\|_{2}^{2}.
\end{align*}

\parindent=8mm \vspace{.1in}
We now derive the logarithmic uncertainty principle associated with the  quaternion shearlet transform.

\parindent=0mm \vspace{.1in}

{\bf Theorem 4.2.}  {\it Let $ F,\Psi$ belong to the  Schwartz class $S(\mathbb R^{2},\mathbb H)$, then the following logarithmic estimate of the uncertainty inequality hold:}
\begin{align*}
\int_{\mathbb R^{2}} \int_{\mathbb R}\int_{\mathbb R^{+}}\ln\|t\|\,\|{\mathcal S_{\Psi}}F(a,s,t)\|_{\mathbb H}^{2}\,\dfrac{da\,ds\,dt}{a^{3}}+ C_{\Psi}\int_{\mathbb R^{2}}\|\omega\|\,{\|{\hat F}(\omega)\|_{\mathbb H}^{2}}\,d\omega\\\geq\,C_{\Psi}\left(\dfrac{\Gamma^{\prime}(1/2)}{\Gamma(1/2)}-\ln\pi\right){{\|F\|}_{2}^{2}}
\end{align*}

\parindent=0mm \vspace{.0in}
 {\it Proof.}   For any non-zero  function  $ F\in S(\mathbb R^{2},\mathbb H)$, the time and frequency spreads satisfy the following inequality  {\cite B} 
 \begin{align*}
\int_{\mathbb R^{2}}\ln\|t\|\,\|F(t)\|^{2}dt+\int_{\mathbb R^{2}}\ln\|\omega\|\,{\|{\hat F}(\omega)\|_{\mathbb H}^{2}}\,d\omega\geq\,\left(\dfrac{\Gamma^{\prime}(1/2)}{\Gamma(1/2)}-\ln\pi\right)\int_{\mathbb R^{2}}\|F(t)\|_{\mathbb H}^{2}\,dt.
\end{align*}
Replacing $F(t)$ by  ${\mathcal S_{\Psi}}F(a,s,t)$ in the above inequality, we obtain
 \begin{align*}
\int_{\mathbb R^{2}}\ln\|t\|\,\|{\mathcal S_{\Psi}}F(a,s,t)\|_{\mathbb H}^{2}dt+\int_{\mathbb R^{2}}\ln\|\omega\|{\left\|\mathscr F_{Q}\left[\mathcal S_{\Psi}F(a,s,t)\right]\right\|_{\mathbb H}^{2}}d\omega\\\geq\,\left(\dfrac{\Gamma^{\prime}(1/2)}{\Gamma(1/2)}-\ln\pi\right)\int_{\mathbb R^{2}}\|{\mathcal S_{\Psi}}F(a,s,t)\|_{\mathbb H}^{2}\,dt
\end{align*}
Integrate the above inequality with respect to measure $da\,ds/a^{3}$ and then apply  Fubini theorem, we obtain
\begin{align*}
\int_{\mathbb R^{2}}\int_{\mathbb R}\int_{\mathbb R^{+}}\ln\|t\|\,\|{\mathcal S_{\Psi}}F(a,s,t)\|_{\mathbb H}^{2}\,\dfrac{da\,ds\,dt}{a^{3}}+\int_{\mathbb R^{2}}\int_{\mathbb R}\int_{\mathbb R^{+}}\ln\|\omega\|{\left\|\mathscr F_{Q}\left[\mathcal S_{\Psi}F(a,s,t)\right]\right\|_{\mathbb H}^{2}}\,\dfrac{da\,ds\,d\omega}{a^{3}}\\\geq\,\left(\dfrac{\Gamma^{\prime}(1/2)}{\Gamma(1/2)}-\ln\pi\right)\int_{\mathbb R^{2}}\int_{\mathbb R}\int_{\mathbb R^{+}}\|{\mathcal S_{\Psi}}F(a,s,t)\|_{\mathbb H}^{2}\,\dfrac{da\,ds\,dt}{a^{3}}.\tag{4.3}
\end{align*}
By estimating the second integral in the L.H.S of the above inequality, we have
\begin{align*}
&\int_{\mathbb R^{2}}\int_{\mathbb R}\int_{\mathbb R^{+}}\ln\|\omega\|{\left\|\mathscr F_{Q}\left[\mathcal S_{\Psi}F(a,s,t)\right]\right\|_{\mathbb H}^{2}}\dfrac{da\,ds\,d\omega}{a^{3}}\\&=\int_{\mathbb R^{2}} \int_{\mathbb R}\int_{\mathbb R^{+}}\ln\|\omega\|\,\,\|\hat{F}(\omega)\|_{\mathbb H}^{2}\,\|\hat{\Psi}(\omega S_{s} A_{a})\|_{\mathbb H}^{2}\,\dfrac{da\,ds\,d\omega}{a^{3/2}}\\&=\int_{\mathbb R^{2}}\left\{\int_{\mathbb R}\int_{\mathbb R^{+}}\dfrac{\|\hat{\Psi}(\omega S_{s} A_{a})\|_{\mathbb H}^{2}}{a^{3/2}}\,da\,ds\right\}\ln\|\omega\|\,\,\|\hat{F}(\omega)\|_{\mathbb H}^{2}\,d\omega\\&=C_{\Psi}\int_{\mathbb R^{2}}\ln\|\omega\|\,\,\|\hat{F}(\omega)\|_{\mathbb H}^{2}\,d\omega.
\end{align*}
Plugging the above estimate in (4.3) , we get
\begin{align*}
\int_{\mathbb R^{2}}\int_{\mathbb R}\int_{\mathbb R^{+}}\ln\|t\|\,\|{\mathcal S_{\Psi}}F(a,s,t)\|_{\mathbb H}^{2}\,\dfrac{da\,ds\,dt}{a^{3}}+C_{\Psi}\int_{\mathbb R^{2}}\ln\|\omega\|\,\,\|\hat{F}(\omega)\|_{\mathbb H}^{2}\,d\omega\geq\,\\\left(\dfrac{\Gamma^{\prime}(1/2)}{\Gamma(1/2)}-\ln\pi\right)\int_{\mathbb R^{2}}\int_{\mathbb R}\int_{\mathbb R^{+}}\|{\mathcal S_{\Psi}}F(a,s,t)\|_{\mathbb H}^{2}\,\dfrac{da\,ds\,dt}{a^{3}}.
\end{align*}
Finally, using the Corollary 3.7,  we get the desired result 
\begin{align*}
\int_{\mathbb R^{2}} \int_{\mathbb R}\int_{\mathbb R^{+}}\ln\|t\|\,\|{\mathcal S_{\Psi}}F(a,s,t)\|_{\mathbb H}^{2}\,\dfrac{da\,ds\,dt}{a^{3}}+C_{\Psi}\int_{\mathbb R^{2}}\|\omega\|\,{\|{\hat F}(\omega)\|_{\mathbb H}^{2}}\,d\omega\\\geq\,C_{\Psi}\left(\dfrac{\Gamma^{\prime}(1/2)}{\Gamma(1/2)}-\ln\pi\right){{\|F\|}_{2}^{2}}.
\end{align*}

\parindent=0mm \vspace{.1in}
{\bf{References}}


\begin{enumerate}

{\small{
\bibitem{D}  Debnath, L., Shah, F.A.:  {\it Wavelet Transforms and Their Applications,} Birkh\"{a}user, New York, (2015).
\bibitem{M}  Daubeachies, I.: {\it Ten Lectures on Wavelets}, SIAM, Philadelphia, (1992).
\bibitem{G}  Kutyniok, G.,  Labate, D: Resolution of the wavefront set using shearlets, preprint (2005).
\bibitem{C}  Chen,L., Kou, K., Liu, M: {\it J. Math. Anal. Appl.} 423(2015) 681-700 (2014).
\bibitem{F}  Folland, G.B., Sitaram, A.:  {\it J. Fourier Anal. Appl.} 3, 207-238 (1997).
\bibitem{B}  Beckner, W.:  {\it Proc. Amer. Math. Soc.} 123, 1897-1905 (1995).
\bibitem{S}  Y. Su, {\it J. Nonlinear Sci. Appl.} 9, 778-786 (2016).

\bibitem{H}  He, J., Yu, B.: {\it Appl. Math. Lett.} 17, 111-121 (2004).
   }}

\end{enumerate}
\end{document}